\documentclass[11pt]{article}
\usepackage{amsmath,amssymb,latexsym, color}
\usepackage{theorem}
\newtheorem{theorem}{Theorem}
\newtheorem{lemma}{Lemma}
\newtheorem{corollary}{Corollary}
\newtheorem{proposition}{Proposition}
{\theorembodyfont{\rmfamily} }

\newcommand{\supp}{\operatorname{supp}}

\title{Connection formulas for general discrete Sobolev polynomials.\
Mehler-Heine  asymptotics}
\author{A. Pe\~{n}a   \thanks{Both authors partially supported by Ministerio de Econom\'{\i}a y Competitividad of Spain under Grant
MTM2012--36732--C03--02 and Diputaci\'on General de Arag\'on project E--64 
(Spain)}\,\,, M. L. Rezola$^{*}$
 \\ \\ {\small Departamento de Matem\'{a}ticas and IUMA.}\\ {Universidad
de Zaragoza (Spain).}}

\date{}

\begin{document}

\maketitle

\begin{abstract} In this paper the discrete Sobolev inner product $$\langle p,q \rangle =\int  p(x) q(x) \,d\mu + \sum_{i=0}^r M_i \, p^{(i)}(c) \, q^{(i)}(c)$$  is considered, where $\mu$ is a finite positive Borel measure supported on an infinite subset of the real line,  $c\in\mathbb{R}$ and $\, M_i \ge 0, \, i = 
0, 1,  \dots, r.$

Connection formulas  for the orthonormal polynomials associated with $\langle . , . \rangle$ are obtained.  As a consequence,  for a wide class of measures $\mu$, we give the  Mehler-Heine asymptotics in the case of the point $c$ is a hard edge of the support of $\mu$. In particular, the case of a symmetric measure $\mu$ is analyzed. Finally, some examples are presented.
\end{abstract}

2000MSC: 42C05; 33C45

\indent Heine type formulas; Bessel functions; Asymptotic zero distribution.

\bigskip

Corresponding author: Ana Pe\~{n}a

Departamento de Matem\'aticas. Universidad de Zaragoza

50009-Zaragoza (Spain)

e-mail: anap@unizar.es

Phone: 34876553224, Fax: 34 976761338

\newpage

\section {Introduction}

Let  $\{p_n(x) \}_{n\ge0}$ be the sequence of orthonormal polynomials with respect to a finite positive Borel measure $\mu$, supported on an infinite subset of $\mathbb{R}$. We denote by  $\{q_n(x) \}_{n\ge0}$ the sequence of orthonormal polynomials with respect to an inner product of the form
\begin{equation}\label{innerproduct}
\langle p,q \rangle =\int  p(x)q(x) \,d\mu + \sum_{i=0}^r M_i \, p^{(i)}(c)  \, q^{(i)}(c),  \quad c \in\mathbb{R}
\end{equation}
where  $\, M_i \ge 0, \, i= 0, \dots, r.$  Such  inner products are called  discrete  Sobolev or Sobolev type and  they have been considered in different contexts. 

In this paper we focus our attention on Mehler-Heine asymptotics for discrete Sobolev orthogonal polynomials.  This asymptotic give us one of the main differences that can be established in order to show how the addition of the derivatives in the inner product influences the orthogonal system.  The Mehler-Heine formulas describe the asymptotic behavior of orthogonal polynomials near the hard edge, i.e. those endpoints of the support of the zero distribution which are also endpoints of the support of the measure. Thus, our interest is to show how the presence  of the masses in the inner product changes the asymptotic behavior around this point.

To prove the Mehler-Heine formula for  Jacobi and Laguerre polynomials, usually the explicit representation of these polynomials are used (see \cite{sz}). Although the situation is a bit different, we would like to mention that recently, in \cite{W},  for some classical multiple orthogonal polynomials, asymptotic formulas of Mehler-Heine type are obtained  using the explicit expression of the polynomials and the Lebesgue's dominated convergence theorem. But there are many others polynomials for which we have not an explicit representation. For example,  Aptekarev in \cite{A} realizes that for certain classes of weight functions supported on $[-1,1]$, the  Mehler-Heine asymptotic  formula depends on the local behavior at the endpoint of the interval of orthogonality. So, this formula has been extended to a  broader class of measures belonging to the Nevai's class.  This result has been applied to deduce the Mehler-Heine formula for  the generalized Jacobi polynomials (see \cite{F}).  On the other hand, for exponential weights  (see \cite{ampr}),  it has been proved the Mehler-Heine formula for the so-called Freud polynomials using the asymptotic formula given by Kriecherbauer and McLaughlin in \cite{K-Mc}  obtained by the Riemann-Hilbert method.

For discrete Sobolev orthogonal polynomials, it is difficult to apply the same method as Jacobi  and Laguerre. This analytic idea was developed for the discrete Laguerre Sobolev orthogonal polynomials in \cite{ampr2011}. There, the authors obtained a new and specific formula for the derivatives of $q_n$ which leads to achieve a uniform  bound in order to use the Lebesgue's dominated convergence theorem. But in a general case, this is quite complicated.

On the other hand, an important tool to get asymptotics is the knowledge of certain connection formulas for $q_n$ in terms of standard polynomials related with $p_n$.  

One of them can be deduced from the well known fact that $\{q_n(x) \}_{n\ge0}$ is quasi--orthogonal of order $r+1$  and consequently  we can  express $q_n$   as a linear combination of the standard orthogonal polynomials  corresponding to the modified measure  $(x-c)^{r+1}\,d\mu(x)$. This connection formula has proved to be fruitful, for example, in the study of relative asymptotics when $\mu$ has compact support, see  \cite{LMV} and  \cite{rms} in a more general setting.  However,  the situation is quite different in the case of measures with unbounded support. So, in \cite{ampr2011},  it was shown that this connection formula is not the adequate to study neither relative asymptotics nor Mehler-Heine formula when $\mu$ is the Laguerre weight.

For  discrete Laguerre Sobolev orthogonal polynomials,  another connection formula was given by Koekoek in  \cite{k1}  with  an arbitrary $r$ and  $c=0$.  This formula has turned out to be of great importance to generating Mehler-Heine asymptotics.  So, in \cite{am},  using the explicit expression of the connection coefficients given in \cite{km} for $r=1$, the authors prove the Mehler-Heine asymptotic for the corresponding Laguerre Sobolev polynomials.  This idea has also been used in \cite{mzfh} for $r=1$ and $c<0$.   However, in an inner product with an arbitrary (finite) number of terms in the discrete part, the problem is that we have not the explicit expression of the coefficients. In spite of this, in \cite{pr} the authors achieve the Mehler-Heine formula for  an arbitrary number of masses and $c=0$.

 In this paper we prove that for a wide class of  measures with support bounded or not, an arbitrary number of masses in the inner product (\ref{innerproduct})  and without taking into account the location of the point $c$ with respect to the support of $\mu$, there exists a  connection  formula for $q_n(x)$  in terms of some canonical transformation of  the  polynomials $p_n$, called Christoffel perturbations.  More precisely,
$q_n(x)=\sum_{j=0}^{r+1} \lambda_{j,n}(x-c)^j\,p_{n-j}^{[2j]}(x)$  were  $\{p_n^{[j]}(x) \}_{n\ge0}$ denotes the sequence of orthonormal polynomials with respect to the measure $\mu_{j}$ with $d\mu_{j}(x)=(x-c)^{j}\,d\mu(x)$.  The main contribution is that we are able to give information  of  asymptotic behavior of the connection coefficients, without the explicit expression of them. This is a significant improvement compared with the previous works. Our interest is focussed in the application of this  connection formula for obtaining the  Mehler-Heine asymptotics and so to prove that whenever the asymptotic behavior near the hard edge involves Bessel functions, the presence of positive masses in the inner product produces a convergence acceleration to this point of $r+1$ zeros of the Sobolev polynomials.

 Finally, we would like to notice that the connection formula obtained has interest by itself because it may be used to get another results as Cohen type inequality  and other asymptotic properties.

The structure of the paper is as follows. In Section 2,  we establish  a connection formula for  orthonormal polynomials with respect to the inner product  (\ref{innerproduct}), where $c$ is an arbitrary real number. Moreover, we show a technical lemma that besides giving us asymptotic behavior at the point  $c$ of the successive derivatives of the polynomials $q_n$ and so of  their kernels, also provides  information about the asymptotic behavior of the  coefficients in  the connection formula.  In Section 3, for a wide class of measures $\mu$, we obtain Mehler-Heine asymptotics for the sequence  $\{q_n(x) \}_{n\ge0}$ where in (\ref{innerproduct}) all the masses are positive and the point $c$ is either an endpoint of the interval where the measure $\mu$ is supported or the origin if the measure is symmetric. As an application, we obtain an important information about the distribution of the zeros of the polynomials $q_n$. In the last section we present some examples to illustrate the theory given.

Throughout this paper we use the notation $x_n \cong y_n$ when the sequence $x_n/y_n$ converges to $1.$ 

\section {Connection formulas}

 Let $\{p_n(x) \}_{n\ge0}$ be the sequence of orthonormal polynomials with respect to the measure  $\mu$ and  $\{q_n(x) \}_{n\ge0}$ the sequence of orthonormal polynomials with respect to the inner product (\ref{innerproduct}). 
 
 In this section we will establish a connection formula for the discrete Sobolev polynomials  $q_n$ in terms of the polynomials $p_n^{[2j]}(x)$ orthogonal with respect to the measure 
$d\mu_{2j}(x)=(x-c)^{2j} \,d\mu(x)$.

\begin{theorem}\label{koekoek}  Assume that the polynomials  $\{p_n(x) \}_{n\ge0}$ satisfy
\begin{equation}\label{condicion-inicial}
p_n(c)\,p_{n-1}^{[2]}(c)\, \dots \,p_{n-(r+1)}^{[2(r+1)]}(c) \not=0
\end{equation}
then there exists  a family of  coefficients  $( \lambda_{j,n})_{j=0}^{r+1}$ , not identically zero,  such that  the following connection formula holds
\begin{equation}\label{connection-formula}
q_n(x)=\sum_{j=0}^{r+1} \lambda_{j,n}\,(x-c)^j\,p_{n-j}^{[2j]}(x), \quad n\ge r+1.
\end{equation}
\end{theorem}

\textbf{Proof.} We will show that there exists  a family of  coefficients  $( \lambda_{j,n})_{j=0}^{r+1},$  not identically zero,  such that  the polynomial $r_n(x)$ defined by  $r_n(x)=\sum_{j=0}^{r+1} \,\lambda_{j,n}\,(x-c)^j\,p_{n-j}^{[2j]}(x)$ satisfies
 \begin{equation}\label{ortogonalidad r_n}
 \langle r_n(x), (x-c)^k\rangle=0, \quad 0\le k \le n-1.
 \end{equation}
  Indeed, for $0 \le j \le r+1 \le k \le n-1$
\begin{align*}
\langle (x-c)^j\,p_{n-j}^{[2j]}(x), (x-c)^k\rangle&=\int(x-c)^jp_{n-j}^{[2j]}(x)\,(x-c)^k\,d\mu(x)\,\\
&=\int\,(x-c)^{k-j}p_{n-j}^{[2j]}(x)\,\,d\mu_{2j}(x)=0.
\end{align*}
Thus, if $0\le k \le r,$  (\ref{ortogonalidad r_n}) leads to the following system 
$$\sum_{j=0}^{r+1} \,\lambda_{j,n}\,\langle (x-c)^j\,p_{n-j}^{[2j]}(x), (x-c)^k\rangle=0$$
 of $r+1$ equations on $r+2$ unknowns  and then we can affirm that it has a non trivial solution $( \lambda_{j,n})_{j=0}^{r+1}$.

 To assure that $r_n(x)=q_n(x)$ it is enough  to prove that  the polynomial $r_n$ has  degree $n$.
Indeed,  if  $\text{deg} \, r_n < n$, the condition (\ref{ortogonalidad r_n}) yields  $r_n\equiv 0$, but this is in contradiction with the hypothesis  (\ref{condicion-inicial}), because if we denote by $\lambda_{j_0,n}$  the first coefficient non zero, then 
$r_n^{(j_0)}(c)=\lambda_{j_0,n}\,j_0! \, p_{n-j_0}^{[2j_0]}(c)\not=0. $    $\quad\Box$

\medskip

\noindent \textbf{Remark 1} If $c$ is not in the interior of the convex hull of the support of the measure $\mu$, the condition (\ref{condicion-inicial}) is always true.

\bigskip

In the next lemma, we get asymptotic estimates of the derivatives $q_n^{(k)}(c)$ from the ones of $p_n^{(k)} (c)$,  which will play  an important role along this paper.

\begin{lemma}\label{cocientederivadas}
 Suppose that  there exists  a strictly increasing function $f$ with $2f(0)+1>0$ and  such that the polynomials $\{p_n(x) \}$ satisfy  the condition
\begin{equation}\label{hipotesis}
p_n^{(k)}(c)\cong C_k(-1)^{n}\, n^{f(k)}, \quad 0\le k \le n.
\end{equation}

Then the following statement holds:
\begin{equation}\label{cociente derivadas}
 \frac{q_{n}^{(k)}(c)}{p_n^{(k)}(c)}\cong \left\{
                              \begin{array}{ll}
                                \displaystyle{\frac{C_k}{n^{2f(k)+1}}}, & \hbox{for k such that $ M_k>0$;} \\
                                \quad C_k, & \hbox{otherwise,}
                              \end{array}
                            \right.
\end{equation}
where $C_k$ is a nonzero constant independent of $n$,  but possibly different in each occurrence.

\end{lemma}

\textbf{Proof.}  We will prove the result by  an induction process concerning the number of positive masses in the inner product (\ref{innerproduct}).

We take the first mass which is positive, namely $M_{j_1}$ ($j_1 \ge 0$),  and consider the sequence of  orthonormal polynomials $\{q_{n,1}\}_{n \ge0}$ with respect to the inner product

\begin{equation*}
(p,q )_{1} = \int p(x)q(x) \,d\mu(x)+  M_{j_1} \, p^{(j_1)}(c) \, q^{(j_1)}(c),
\end{equation*}
where  $q_{n,1}(x)=\tilde{\gamma}_{n,1}\,x^n+\dots$.

The Fourier expansion of the polynomial $q_{n,1}$ in the orthonormal basis $\{p_n\}_{n\ge0}$ ($p_n(x)=\gamma_n\,x^n+\dots $)  leads to
\begin{equation*}
q_{n,1}(x) = \frac{\tilde{\gamma}_{n,1}}{\gamma_n}p_n(x) - M_{j_1}\,  q_{n,1}^{(j_1)}(c)\,  K_{n-1}^{(0,j_1)}(x,c) \,,
\end{equation*}
and therefore
\begin{equation}\label{qn,1(x)}
q_{n,1}(x) = \frac{\tilde{\gamma}_{n,1}}{\gamma_n}\left[p_n(x) - \frac{M_{j_1} \, p_n^{(j_1)}(c)}{1+ M_{j_1}   K_{n-1}^{(j_1,j_1)}(c,c)} K_{n-1}^{(0,j_1)}(x,c)\right] \,
\end{equation}
where, as usual, we denote by $K_n^{(k,h)}(x,y)$ the derivatives of the $n$th kernel for the sequence
$\{p_n\}_{n\ge0}$ 
 $$K_n^{(k,h)}(x,y)= \frac{\partial^{k+h}}
{\partial x^k \, \partial x^h } K_n (x,y) = \sum _{i=0}^n p_i^{(k)}(x) \, p_i^{(h)}(y), \quad k,h  \in \mathbb{N}\cup \{0\}.$$

On the other hand,  for $0\le k\le n$, the hypothesis  for the function $f$ allows us to affirm that $ f(k)+f(j_1)+1>0$. So,   applying Stolz criterion (see, e.g \cite{k}) and  the hypothesis (\ref{hipotesis}) we obtain
 \begin{equation}\label{Kn(k,j1)}
\lim_n \frac{K_{n}^{(k,j_1)}(c,c)}{n^{f(k)+f(j_1)+1}}=\lim_n \frac{p_n^{(k)}(c)\, p_n^{(j_1)}(c)}{(f(k)+f(j_1)+1)\,n^{f(k)+f(j_1)}}=\frac{C_{k}C_{j_1}}{f(k)+f(j_1)+1} \not=0,
\end{equation}
and so 
 \begin{equation}\label{Kn(k,j1)1}
\frac{K_{n}^{(k,j_1)}(c,c)}{n\, p_n^{(k)}(c)\, p_n^{(j_1)}(c)} \cong  \frac{1}{f(k)+f(j_1)+1}
\end{equation}
Moreover, it is easy to check  that
$$\left( \frac{\gamma_n}{\tilde{\gamma}_{n,1}} \right)^2=\left[1+\frac{M_{j_1}(p_n^{(j_1)}(c))^2}{1+ M_{j_1} K_{n-1}^{(j_1,j_1)}(c,c)}\right]$$ and thus
\begin{equation}\label{equivnormas1}
\frac{\tilde{\gamma}_{n,1}}{\gamma_n}\cong 1.
\end{equation}

Now, taking derivatives $k$ times in (\ref{qn,1(x)}) and evaluating at $x=c$, we have
\begin{equation*}\label{qn,1/pn}
\frac{q_{n,1}^{(k)}(c)}{p_n^{(k)}(c)} =\frac{\tilde{\gamma}_{n,1}}{\gamma_n}\left[1- \frac{M_{j_1}K_{n-1}^{(k,j_1)}(c,c)}{1+M_{j_1}K_{n-1}^{(j_1,j_1)}(c,c)}\frac{p_n^{(j_1)}(c)}{p_n^{(k)}(c)} \right] \,.
\end{equation*}
Then, by (\ref{Kn(k,j1)})  and  (\ref{equivnormas1}) , we get 

\begin{equation*}
\frac{q_{n,1}^{(j_1)}(c)}{p_n^{(j_1)}(c)}=\frac{\tilde{\gamma}_{n,1}}{\gamma_n}\,\frac{1}{1+M_{j_1}K_{n-1}^{(j_1,j_1)}(c,c)}\cong \frac{C_{j_1}}{n^{2f(j_1)+1}},
\end{equation*}
and for $k\not=j_1$, from (\ref{Kn(k,j1)1})  and  (\ref{equivnormas1})
\begin{equation*}
\frac{q_{n,1}^{(k)}(c)}{p_n^{(k)}(c)}\cong \left[1-\frac{2f(j_1)+1}{f(k)+f(j_1)+1} \right]=\frac{f(k)-f(j_1)}{f(k)+f(j_1)+1}\not=0.
\end{equation*}

 If there are no more positive masses, since $q_{n,1}=q_n$  we have concluded the proof. Otherwise, suppose that the result  holds for the sequence of orthonormal polynomials $\{q_{n,s-1}\}_{n \ge0}$ orthogonal with respect to the inner product
\begin{align*}
(p,q )_{s-1} &= \int  p(x)q(x) \,d\mu(x) + \sum_{i=1}^{s-1} M_{j_i}  \, p^{(j_{i})}(c)  \, q^{(j_{i})}(c),
\end{align*}
where $j_1 <j_2 < \dots <j_{s-1}$ and all these masses are positive.

Now, we have to prove the result for the orthonormal polynomials $q_{n,s}$  with respect to
\begin{equation*}
(p,q )_{s} =(p,q )_{s-1} +M_{j_{s}} \, p^{(j_{s})}(c)  \, q^{(j_{s})}(c),
\end{equation*}
where $M_{j_{s}}>0$, and we can work as before. Then the Fourier expansion of the polynomial $q_{n,s}$  ($q_{n,s}(x)=\tilde{\gamma}_{n,s}\,x^n+\dots $)  in the orthonormal basis $\{q_{n,s-1}\}_{n\ge0}$ leads to
\begin{equation*}
q_{n,s}(x) = \frac{\tilde{\gamma}_{n,s}}{\tilde{\gamma}_{n,s-1}}\, q_{n,s-1}(x) - M_{j_s}\,  q_{n,s}^{(j_s)}(c) \, K_{n-1, s-1}^{(0,j_s)}(x,c) \,,
\end{equation*}
where  $K_{n,s-1}$ denotes the corresponding $n$th kernel for the sequence $\{q_{n,s-1}\}$ and $$K_{n,s-1}^{(k,h)}(x,y)=
\displaystyle \sum _{i=0}^n
q_{i,s-1}^{(k)}(x) \, q_{i,s-1}^{(h)}(y), \quad k,h  \in \mathbb{N}\cup \{0\}.$$
Therefore, 
\begin{equation}\label{qn,s(x)}
q_{n,s}(x) =  \frac{\tilde{\gamma}_{n,s}}{\tilde{\gamma}_{n,s-1}}\left[ q_{n,s-1}(x) - \frac{M_{j_s} \, q_{n,s-1}^{(j_s)}(c)}{1+ M_{j_s}  K_{n-1,s-1}^{(j_s,j_s)}(c,c)} K_{n-1,s-1}^{(0,j_s)}(x,c)  \right]\,,
\end{equation}
and
\begin{equation}\label{cociente gammas}
\left( \frac{\tilde{\gamma}_{n,s-1}}{\tilde{\gamma}_{n,s}} \right)^2=\left[1+\frac{M_{j_s} \, (q_{n,s-1}^{(j_s)}(c))^2}{1+ M_{j_s}  K_{n-1,s-1}^{(j_s,j_s)}(c,c)}\right].
\end{equation}

Applying Stolz criterion, the hypothesis for the function $f$  and
 the hypothesis  for $\{ q_{n,s-1} \}_{n\ge0}$, we can obtain
\begin{equation*}\label{K_{n-1,s-1}}
 K_{n,s-1}^{(k,j_s)}(c,c) \cong \left\{
                              \begin{array}{ll}
                                C_k \, n^{f(k)+f(j_s)+1}, & \hbox{if $k \not=j_1,\dots, j_{s-1}$;} \\
                                C_k \, n^{f(j_s)-f(k)}, & \hbox{if $ k=j_1,\dots, j_{s-1}$,}
                              \end{array}
                            \right.
\end{equation*}
where $C_k$ is a nonzero constant independent of $n$,  but possibly different in each occurrence.

Indeed, for $k \not=j_1,\dots, j_{s-1}$,
\begin{align}\label{K_{n,s-1}^{(k,js}(c,c)}
&  \lim_n \frac{K_{n,s-1}^{(k,j_s)}(c,c)}{n^{f(k)+f(j_s)+1}}
 =\lim_n \frac{q_{n,s-1}^{(k)}(c)\,q_{n,s-1}^{(j_s)}(c)}{(f(k)+f(j_s)+1)\,n^{f(k)+f(j_s)}} \\ \nonumber
&=\lim_n \frac{q_{n,s-1}^{(k)}(c)}{p_{n}^{(k)}(c)}\lim_n\frac{q_{n,s-1}^{(j_s)}(c)}{p_n^{(j_s)}(c)}\lim_n  \frac{p_n^{(k)}(c) \, p_n^{(j_s)}(c)}{(f(k)+f(j_s)+1)\, n^{f(k)+f(j_s)}}\not=0,  \\ \nonumber
\end{align}
and, for $k=j_1,\dots, j_{s-1}$,
\begin{align}\label{K_{n,s-1}^{(k,js)}(c,c)}
& \lim_n \frac{K_{n,s-1}^{(k,j_s)}(c,c)}{n^{f(j_s)-f(k)}}
=\lim_n  \frac{q_{n,s-1}^{(k)}(c) \, q_{n,s-1}^{(j_s)}(c)}{(f(j_s)-f(k)) \,\, n^{f(j_s)-f(k)-1}} \\ \nonumber
&=\lim_n \frac{p_n^{(k)}(c) \, p_n^{(j_s)}(c)} {(f(j_s)-f(k)) \,n^{f(k)+f(j_s)}} \lim_n  n^{2f(k)+1} \frac{q_{n,s-1}^{(k)}(c)}{p_n^{(k)}(c)} \lim_n \frac{q_{n,s-1}^{(j_s)}(c)}{p_{n}^{(j_s)}(c)}\not=0.
\end{align}

Then from (\ref{cociente gammas})
\begin{equation}\label{equivnormas s}
\frac{\tilde{\gamma}_{n,s}  }{\tilde{\gamma}_{n,s-1}  }\cong 1.
\end{equation}

Now, taking derivatives $k$ times in (\ref{qn,s(x)})  and evaluating at $x=c$, we obtain
\begin{equation}\label{qn,s/pn}
\frac{q_{n,s}^{(k)}(c)}{p_n^{(k)}(c)}\cong \,\frac{q_{n,s-1}^{(k)}(c)}{p_n^{(k)}(c)}\left[1 - \frac{M_{j_s}K_{n-1,s-1}^{(k,j_s)}(c,c)}{1+M_{j_s}K_{n-1,s-1}^{(j_s,j_s)}(c,c)}
\frac{q_{n,s-1}^{(j_s)}(c)}{q_{n,s-1}^{(k)}(c)} \right] \,.
\end{equation}

For $k=j_s$,  the hypothesis for $q_{n,s-1}$ and the estimation of the kernel yield 

\begin{equation*}
\frac{q_{n,s}^{(j_s)}(c)}{p_n^{(j_s)}(c)}\cong \,\frac{q_{n,s-1}^{(j_s)}(c)}{p_n^{(j_s)}(c)}
\frac{1}{1+M_{j_s}K_{n-1,s-1}^{(j_s,j_s)}(c,c)}\cong \frac{C_{j_s}}{n^{2f(j_s)+1}},
\end{equation*}
with $C_{j_s}$ a nonzero constant independent of $n$.

 Moreover, for $k\not=j_s$, taking into account (\ref{K_{n,s-1}^{(k,js}(c,c)}), (\ref{K_{n,s-1}^{(k,js)}(c,c)}) and the hypothesis for $q_{n,s-1}$, we can deduce

\begin{equation}\label{cociente nucleos}
\frac{K_{n-1,s-1}^{(k,j_s)}(c,c)}{K_{n-1,s-1}^{(j_s,j_s)}(c,c)}\frac{q_{n,s-1}^{(j_s)}(c)}{q_{n,s-1}^{(k)}(c)}
 \cong \left\{
\begin{array}{ll}
 \frac{2f(j_s)+1}{f(k)+f(j_s)+1}, & \hbox{if $k\not=j_1,\dots, j_{s-1}$;} \\
 \frac{2f(j_s)+1}{f(j_s)-f(k)}, & \hbox{if $k=j_1,\dots, j_{s-1}$.}
\end{array}
\right.
\end{equation}

Thus, taking limits in (\ref{qn,s/pn}), we get for the polynomials $q_{n,s}$,
\begin{equation*}
 \frac{q_{n,s}^{(k)}(c)}{p_n^{(k)}(c)}\cong \left\{
                              \begin{array}{ll}
                                \displaystyle{\frac{C_k}{n^{2f(k)+1}}}, & \hbox{if $ k=j_1,\dots, j_{s}$;} \\
                                \quad C_k, & \hbox{otherwise,}
                              \end{array}
                            \right.
\end{equation*}
where the hypothesis for $f$ allows us to affirm that $C_k$ is a nonzero constant independent of $n$,  but possibly different in each occurrence. Hence the result follows.$\quad \Box$

\bigskip

\begin{corollary}\label{cocientederivadas r+1}  Under the same hypothesis of the previous lemma we have
\begin{equation}
\frac{q_{n}^{(r+1)}(c)}{p_n^{(r+1)}(c)}\cong \prod_{i=1}^s \frac{f(r+1)-f(j_i)}{f(r+1)+f(j_i)+1} 
\end{equation}
where  $(M_{j_i})_{i=1}^s$ are the positive masses in the inner product (\ref{innerproduct}).
\end{corollary}

\textbf{Proof.}  We obtain the result applying a recursive process concerning the number of positive masses in the inner product (\ref{innerproduct}).

Indeed, using (\ref{qn,s/pn}) and  (\ref{cociente nucleos}), it follows that
$$\frac{q_{n,s}^{(r+1)}(c)}{p_{n}^{(r+1)}(c)} \cong \frac{q_{n,s-1}^{(r+1)}(c)}{p_{n}^{(r+1)}(c)}\left[ \frac{f(r+1)-f(j_s)}{f(r+1)+f(j_s)+1} \right]$$  
and we get the result. $\quad \Box$

\medskip

\noindent \textbf{Remark 2}  Lemma \ref{cocientederivadas} and hence Corollary \ref{cocientederivadas r+1}  are also true if in the condition (\ref{hipotesis}) the factor $(-1)^n$ is deleted.

\bigskip

Next, as a consequence of Lemma \ref{cocientederivadas}, we  prove that  whenever the polynomials $p_n^{[2j]}(x)$ satisfy  a  similar condition to (\ref{hipotesis}), then  there  exists  the connection formula (\ref{connection-formula})  and moreover there exists  limit of their connection coefficients  $\lambda_{j,n}, \,j=0, 1, \dots r+1$.

\begin{theorem} \label{limitecoeficientesKoekoek}   Suppose that there exists  a strictly increasing function  $f$ with $2f(0)+1>0$   and such that for all $j=0, 1,  \dots, r+1$,  the polynomials $\{p_n^{[2j]}(x)\}$  satisfy  the condition
\begin{equation}\label{hipotesis-j}
(p_n^{[2j]})^{(k)}(c)\cong C_{k,j}\,(-1)^{n}\, n^{f(k+j)}, \quad 0 \le k  \le n,
\end{equation} 
where $C_{k,j}$ is a nonzero constant independent of $n$.

Then, there exists
 $$  \lim_n \,\lambda_{j,n}=\lambda_j \in \mathbb{R}, \quad j=0, 1,  \dots, r+1,$$
 where  $\{\lambda_{j,n}\}_{0}^{r+1}$  are the coefficients in the connection formula (\ref{connection-formula}).
 Moreover,  if all the masses in the inner product (\ref{innerproduct}) are positive, we obtain
 \begin{equation*}
 \lim_n \,\lambda_{j,n}=0,\quad j=0, 1,  \dots, r
 \end{equation*}
and
\begin{equation*}
 \lim_n \,\lambda_{r+1,n}\not=0.
\end{equation*}
\end{theorem}

\textbf{Proof.}  Notice that the existence of   the connection formula (\ref{connection-formula})  for n large enough is a straightforward consequence  of 
(\ref{hipotesis-j}). 

So, taking derivatives $k$ times in (\ref{connection-formula})  and evaluating at $x=c$, we deduce
\begin{equation}\label{cocientederivadasortonormales}
\frac{q_n^{(k)}(c)}{p_n^{(k)}(c)}=\sum_{j=0}^{k}\lambda_{j,n} \left(
\begin{array}{c}
 k \\
 j \\
  \end{array}
   \right)
j!\,A_j(k,n), \quad 0 \le k \le r+1,
\end{equation}
where $A_0(k,n)=1$ and
\begin{equation}\label{expresionA}
A_j(k, n)=\displaystyle \frac{(p_{n-j}^{[2j]})^{(k-j)}(c)}{p_n^{(k)}(c)}.
\end{equation}
From condition (\ref{hipotesis-j}),  we can deduce that there exists $\displaystyle \lim_n A_j(k, n)\not=0$. Then,  applying recursively (\ref{cociente derivadas}) and  (\ref {cocientederivadasortonormales}),  we can assure that   there exists $\displaystyle \lim_n \lambda_{j,n}=\lambda_{j}, \, j=0, 1, \dots r+1$. More precisely,  for $k=0$ we have
$$\lim_n \lambda_{0,n}=\lim_n \frac{q_n(c)}{p_n(c)}=\lambda_{0} \, \left\{
\begin{array}{ll}
=0, & \hbox{if $M_0>0$;} \\
\not=0, & \hbox{if $M_0=0$.}
 \end{array}
 \right.$$
\noindent Now, from (\ref{cocientederivadasortonormales}) for $k=1$ and  (\ref{cociente derivadas})  we get

$$ \lim_n \lambda_{1,n}=\lim_n \frac{1}{A_1(1,n)}\left(\frac{q'_n(c)}{p_n'(c)}-\lambda_{0,n}\right)=\lambda_{1}.$$
\noindent Observe that
$$\lambda_{1} \, \left\{
\begin{array}{ll}
=0, & \hbox{if $M_0>0$ and $M_1>0$ ;} \\
\not=0, & \hbox{if $M_0>0$ and $M_1=0$.}
 \end{array}
 \right.$$
 \noindent In this way, recursively, if $M_0 M_1 \dots M_i>0$ and $M_{i+1}=0,$ we can assure that
\begin{equation*}
\lim_n \lambda_{j,n}=\lambda_{j} \, \left\{
\begin{array}{ll}
 =0, & \hbox{if $0 \le j \le i$;} \\
 \not=0, & \hbox{if $j=i+1$,}
 \end{array}\right.
 \end{equation*}
and we obtain the result. $\quad \Box$

\medskip

\noindent \textbf{Remark 3} Theorem \ref{limitecoeficientesKoekoek}   is also true if in the condition (\ref{hipotesis-j})  the factor $(-1)^n$ is deleted.

\bigskip

\subsection {Symmetric case} 

Next, if the measure $ \mu $ is symmetric,  by a symmetrization process, we will obtain similar results to those discussed above.

Let  $\{p_n(x) \}_{n\ge0}$ be the sequence of orthonormal polynomials with respect to a symmetric positive Borel measure $\mu$ on $(-d,d)$ where $0<d\le \infty$.  Let $\nu$ be the image of measure $\mu$ on $J=(0,d^2)$ under the mapping $\Phi(x)=x^2,$ i.e. $\nu=\Phi(\mu)$.

A symmetrization process, see \cite{Ch} for monic polynomials,  yields 
\begin{equation}\label{simetrizacion pn}
p_{2n}(x)=u_n(x^2), \quad p_{2n+1}(x)=x\,u_n^*(x^2),
\end{equation}
where $\{u_n(x) \}_{n\ge0}$ and $\{u_n^*(x) \}_{n\ge0}$ are the sequences of orthonormal polynomials with respect to the measures $d\,\nu(x)$ and $x\,d\nu(x)$, respectively.

Now, we  rename the inner product (\ref{innerproduct}) as 
\begin{equation}\label{innerproduct-simetrico}
\langle p,q \rangle =\int  p(x)q(x) \,d\mu + \sum_{i=0}^{2r+1} M_i \, p^{(i)}(0)  \, q^{(i)}(0), 
\end{equation}
where  $\, M_i \ge 0, \, i = 0, 1, \dots, 2r+1.$ This inner product has been already considered in \cite{ammr}.

First notice  that if the initial measure $\mu$ is symmetric, then the Sobolev type polynomials $q_n$ orthogonal with respect to (\ref{innerproduct-simetrico})  are also symmetric.  Again, by the symmetrization process we can write
\begin{equation}\label{simetrizacion qn}
q_{2n}(x)=s_n(x^2), \quad q_{2n+1}(x)=x\,s_n^*(x^2),
\end{equation}
where now  the sequences of orthonormal polynomials $\{s_n(x) \}_{n\ge0}$ and $\{s_n^*(x) \}_{n\ge0}$ are orthogonal with respect to the Sobolev type inner products: 
\begin{equation}\label{innerproduct-pares}
\langle p,q \rangle_1=\int_J p(x)q(x)\,d\nu(x)+\sum_{i=0}^r  \overline {M_{2i}} \, p^{(i)}(0) \, q^{(i)}(0), 
\end{equation}
and
\begin{equation}\label{innerproduct-impares}
\langle p,q \rangle_2=\int_J p(x)q(x)\,xd\nu(x)+\sum_{i=0}^r \overline {M_{2i+1}} \, p^{(i)}(0) \, q^{(i)}(0),
\end{equation}
respectively, where 
$$\overline {M_{2i}}=\left(\frac{(2i)!}{i!} \right)^2\,M_{2i}, \quad \overline {M_{2i+1}}=\left(\frac{(2i+1)!}{i!} \right)^2\,M_{2i+1}, \quad i=0, \dots, r,$$
(see Theorem 2 in \cite{ammr}).

Also, since the measure $\mu_{2j}$ is symmetric, then the polynomials $p_n^{[2j]}(x)$ are symmetric and it is easy to check that 
 \begin{equation}\label{simetrizacion pn-j}
p_{2n}^{[2j]}(x)=u_n^{[j]}(x^2), \quad p_{2n+1}^{[2j]}(x)=x\,(u_n^*)^{[j]}(x^2),
\end{equation}
where $\{u_n^{[j]}(x) \}_{n\ge0}$ and  $\{(u_n^*)^{[j]}(x) \}_{n\ge0}$ are the sequences of orthonormal polynomials with respect to the measures $x^j\,d\nu(x)$ and $x^{j+1}\,d\nu(x)$, respectively.

\bigskip

\begin{theorem}\label{koekoek-simetrico}  Let $\{p_n(x) \}_{n\ge0}$ be the sequence of orthonormal polynomials with respect to a symmetric measure  $\mu$ and  $\{q_n(x) \}_{n\ge0}$ the sequence of orthonormal polynomials with respect to the inner product (\ref{innerproduct-simetrico}).  

Then, there exist two families of  coefficients  $( \lambda_{j,n})_{j=0}^{r+1}$  and  $( \lambda^*_{j,n})_{j=0}^{r+1}$  not identically zero,  such that  the following connection formulas hold
\begin{equation}\label{connection-formula-pares}
q_{2n}(x)=\sum_{j=0}^{r+1} \lambda_{j,n}\,x^{2j}\,p_{2n-2j}^{[4j]}(x), \quad n\ge r+1,
\end{equation}
and
\begin{equation}\label{connection-formula-impares}
q_{2n+1}(x)=\sum_{j=0}^{r+1} \lambda^*_{j,n}\,x^{2j}\,p_{2n+1-2j}^{[4j]}(x), \quad n\ge r+1.
\end{equation}
\end{theorem}

\textbf{Proof.}  The result  is a simple consequence of Theorem \ref{koekoek} and the symmetrization process described above. Indeed, from (\ref{simetrizacion pn-j})  and since $p_n^{[2j]}(x)$ are symmetric polynomials, 
we get
$$u_n(0)u_{n-1}^{[2]}(0)\dots  u_{n-(r+1)}^{[2(r+1)]}(0)=p_{2n}(0)p_{2n-2}^{[4]}(0) \dots p_{2n-2(r+1)}^{[4(r+1)]}(0) \not=0,$$
and 
\begin{align*}
&u_n^*(0)(u_{n-1}^*)^{[2]}(0)\dots  (u_{n-(r+1)}^*)^{[2(r+1)]}(0)\\ \nonumber
& =(p_{2n+1})'(0)\,(p_{2n-1}^{[4]})'(0) \dots (p_{2n+1-2(r+1)}^{[4(r+1)]})'(0) \not=0.
\end{align*}
Now, taking into account  (\ref{simetrizacion qn})-(\ref{innerproduct-impares}), we can apply Theorem \ref{koekoek} to the polynomials $s_n$ and $s_n^*$ and so there exist two families of  coefficients  $( \lambda_{j,n})_{j=0}^{r+1}$  and  $( \lambda^*_{j,n})_{j=0}^{r+1}$,  not identically zero,  such that 
\begin{equation}\label{connection-formulas-simetrico}
s_n(x)=\sum_{j=0}^{r+1}  \lambda_{j,n}x^j\,u_{n-j}^{[2j]}(x) , \quad \quad  s_n^*(x)=\sum_{j=0}^{r+1}  \lambda^*_{j,n}x^j\,(u_{n-j}^*)^{[2j]}(x). 
\end{equation}
To conclude it is enough to use  again (\ref{simetrizacion qn}) and  (\ref{simetrizacion pn-j}).  $\quad \Box$

\begin{lemma}\label{cocientederivadas-simetrico}  Suppose that there exists  a strictly increasing function  $f$ with $2f(0)+1>0$ and  such that the polynomials $\{p_n\}$ satisfy  the conditions
\begin{equation*}
p_{2n}^{(2k)}(0)\cong C_k(-1)^{n}\, n^{f(2k)}, \quad 0\le k \le n,
\end{equation*}
\begin{equation*}
p_{2n+1}^{(2k+1)}(0)\cong C_k(-1)^{n}\, n^{f(2k+1)}, \quad 0\le k \le n.
\end{equation*}

Then the following statements hold:
\begin{equation*}
 \frac{q_{2n}^{(2k)}(0)}{p_{2n}^{(2k)}(0)}\cong \left\{
                              \begin{array}{ll}
                                \displaystyle{\frac{C_k}{n^{2f(2k)+1}}}, & \hbox{for k such that $ M_{2k}>0$;} \\
                                \quad C_k, & \hbox{otherwise,}
                              \end{array}
                            \right.
\end{equation*}

\begin{equation*}
 \frac{q_{2n+1}^{(2k+1)}(0)}{p_{2n+1}^{(2k+1)}(0)}\cong \left\{
                              \begin{array}{ll}
                                \displaystyle{\frac{C_k}{n^{2f(2k+1)+1}}}, & \hbox{for k such that $ M_{2k+1}>0$;} \\
                                \quad C_k, & \hbox{otherwise,}
                              \end{array}
                            \right.
\end{equation*}
where $C_k$ is a nonzero constant independent of $n$,   but possibly different in each occurrence.
\end{lemma}

\textbf{Proof.} Again we will apply the symmetrization process. 

From (\ref{simetrizacion pn-j}),  with $j=0$,  we have that 
$$u_n^{(k)}(0)=\frac{k!}{(2k)!}\,p_{2n}^{(2k)}(0)\cong C_k (-1)^{n}\, n^{f(2k)}=C_k(-1)^{n}\, n^{g(k)}$$ and 
$$(u_n^*)^{(k)}(0)=\frac{k!}{(2k+1)!}\,p_{2n+1}^{(2k+1)}(0)\cong C_k(-1)^{n}\, n^{f(2k+1)}=C_k(-1)^{n}\, n^{g^*(k)}$$ 
where $g,g^*$ are  strictly increasing functions  satisfying  $2g(0)+1>0,\, 2g^*(0)+1>0.$ Then we can apply Lemma \ref{cocientederivadas} and so we obtain:

\begin{equation*}
 \frac{s_{n}^{(k)}(0)}{u_{n}^{(k)}(0)}\cong \left\{
                              \begin{array}{ll}
                                \displaystyle{\frac{C_k}{n^{2g(k)+1}}}, & \hbox{for k such that $ \overline{M_{2k}}>0$;} \\
                                \quad C_k, & \hbox{otherwise,}
                              \end{array}
                            \right.
\end{equation*}

\begin{equation*}
 \frac{(s_{n}^*)^{(k)}(0)}{(u_{n}^*)^{(k)}(0)}\cong \left\{
                              \begin{array}{ll}
                                \displaystyle{\frac{C_k}{n^{2g^*(k)+1}}}, & \hbox{for k such that $ \overline{M_{2k+1}}>0$;} \\
                                \quad C_k, & \hbox{otherwise,}
                              \end{array}
                            \right.
\end{equation*}

\noindent where $C_k$ is a nonzero constant independent of $n$, but possibly different in each occurrence. 

Thus the result is proved. $\quad \Box$

\begin{theorem} \label{limitecoeficientesKoekoek-simetrico}   Assume that there exists  a strictly increasing function  $f$ with $2f(0)+1>0$  and  such that  for all $j=0, 1,  \dots, r+1$ the polynomials $\{p_n^{[4j]}(x)\}$ satisfy  the conditions
\begin{equation}\label{hipotesis-j-pares}
(p_{2n}^{[4j]})^{(2k)}(0)\cong C_{k,j}\,(-1)^{n}\, n^{f(2k+2j)},  \quad 0\le k \le n,
\end{equation} 

\begin{equation}\label{hipotesis-j-impares}
(p_{2n+1}^{[4j]})^{(2k+1)}(0)\cong C_{k,j}\,(-1)^{n}\, n^{f(2k+2j+1)}, \quad 0\le k \le n,
\end{equation} 
where $C_{k,j}$ is a nonzero constant independent of $n$,   but possibly different in each occurrence.

Then, there exists
 $$  \lim_n \,\lambda_{j,n}=\lambda_j\in \mathbb{R} , \quad \text{and} \quad  \lim_n \,\lambda_{j,n}^*=\lambda_j^*\in \mathbb{R}, \quad \quad j=0, 1, \dots, r+1,$$
 where $\{\lambda_{j,n}\}_{0}^{r+1}, \,\{\lambda_{j,n}^*\}_{0}^{r+1}$  are  the families of coefficients in formulas (\ref{connection-formula-pares}) and (\ref{connection-formula-impares}).
 Moreover,  if all the masses in the inner product (\ref{innerproduct-simetrico}) are positive we obtain
 \begin{equation*}
 \lim_n \,\lambda_{j,n}=\lim_n \,\lambda_{j,n}^*=0,\quad j=0, 1,  \dots, r
 \end{equation*}
and
$$ \lim_n \,\lambda_{r+1,n}\not=0,  \quad \quad  \lim_n \,\lambda_{r+1,n}^*\not=0.$$
\end{theorem}

\textbf{Proof.}   From (\ref{simetrizacion pn-j}) and  (\ref {hipotesis-j-pares})-(\ref {hipotesis-j-impares}) we have 
\begin{align*}
(u_{n}^{[2j]})^{(k)}(0)&=\frac{k!}{(2k)!}(p_{2n}^{[4j]})^{(2k)}(0)\\ \nonumber
&\cong C_{k,j} (-1)^n\, n^{f(2k+2j)}=C_{k,j}(-1)^n\, n^{g(k+j)},
\end{align*}
and 
\begin{align*}
((u_{n}^*)^{[2j]})^{(k)}(0)&=\frac{k!}{(2k+1)!}(p_{2n+1}^{[4j]})^{(2k+1)}(0)\\ \nonumber
&\cong C_{k,j} (-1)^n\, n^{f(2k+2j+1)}=C_{k,j}(-1)^n\, n^{g^*(k+j)},
\end{align*}
where $g,g^*$ are  strictly increasing functions  satisfying  $2g(0)+1>0,\, 2g^*(0)+1>0,$ and  $C_{k,j}$ is a nonzero constant independent of $n$. 

Thus  the result follows, taking into account the connection formulas for  $s_n(x)$ and $s_n^*(x)$ given in  (\ref{connection-formulas-simetrico})
 and  Theorem \ref{limitecoeficientesKoekoek}.  $\quad \Box$

\bigskip

\section {Mehler-Heine type formulas} 

 Let  $\{p_n(x) \}_{n\ge0}$ be the sequence of orthonormal polynomials with respect to a positive measure $\mu$ supported on an interval $I$. Let  $\{q_n(x) \}_{n\ge0}$  be the sequence of Sobolev type orthonormal polynomials with respect to the inner product (\ref{innerproduct}) where all the masses are positive and $c$ is an endpoint of  the interval $I$. We assume without loss of generality that $c=\inf I , \, c\in \mathbb{R}. $
 
In this section we obtain Mehler-Heine asymptotics for  $\{q_n(x) \}_{n\ge0}$.  To do this  we will use essentially the connection formulas  and  the asymptotic estimates of their connection coefficients given in Section 2.

Remind that  the Bessel functions  $J_{\nu}$ of  order  $\nu $, $\nu >-1$ are defined by
\begin{equation*}
J_{\nu}(z)=\sum_{n=0}^{\infty} \frac{(-1)^n}{n!
\,\Gamma(n+\nu+1)} \left( \frac{z}{2} \right)^{2n+\nu}, \quad  z \in \mathbb{C}.
\end{equation*}

\bigskip

\begin{theorem} \label{MH-general}  Suppose that the sequence  $\{p_n^{[2j]}(x) \}_{n\ge0}$ satisfies  uniformly  on  compact  subsets of $\mathbb{C}$, for all $j=0, 1, \dots, r+1$,  the following Mehler-Heine asymptotics

\begin{equation}\label{MH pn-2j}
 \lim_{n}\, (-1)^n\frac{a_n^{1/2}}{b_n^j}\,p_n^{[2j]}(c+\frac{z^2}{b_n})=z^{-(\nu+2j)}\,J_{\nu+2j}(2z),
\end{equation}
 where  
\begin{equation}\label{estimacion an-bn}
 a_n^{-1/2}\cong An^{a}, \quad b_n\cong Bn^b,  \quad A,\,B,\, b >0 , \quad  \nu >-1, 
\end{equation}
and
\begin{equation}\label{ligadura}
2a+1=b(\nu+1).
\end{equation}

Then 
\begin{equation}\label{MH qn}
 \lim_{n}\, (-1)^n\,a_n^{1/2}\,q_n(c+\frac{z^2}{b_n})=(-1)^{r+1}\,z^{-\nu}\,J_{\nu+2r+2}(2z),
 \end{equation}
uniformly  on  compact  subsets of $\mathbb{C}.$
\end{theorem}

\textbf{Proof.}  From hypothesis (\ref{MH pn-2j}) and using the Taylor series for the polynomials  $p_n^{[2j]}(c+\frac{z^2}{b_n})$ at the point $z=c$, we have 
$$ (-1)^n\frac{a_n^{1/2}}{b_n^j}\sum_{k=0}^{n}\frac{(p_n^{[2j]})^{(k)}(c)}{k! \, b_n^k}\,z^{2k} \longrightarrow\sum_{k=0}^{\infty} \frac{(-1)^k}{k! \, \Gamma(k+\nu+2j+1)}\,z^{2k} $$
uniformly  on  compact  subsets of $\mathbb{C},$ which implies 
$$(p_n^{[2j]})^{(k)}(c)\cong (-1)^{n+k}\frac{a_n^{-1/2}b_n^{k+j}}{\Gamma(k+\nu+2j+1)}, \quad  0 \le k \le n.$$
So, writing $f(x)=bx+a$, we have that $f$ is a stricty increasing function with $2f(0)+1=b(\nu+1)>0$  and so the polynomials $p_n^{[2j]}$ satisfy
\begin{equation}\label{ultimo coeficiente}
(p_n^{[2j]})^{(k)}(c)\cong (-1)^{n+k}\,\frac{A\,B^{k+j}}{\Gamma(k+\nu+2j+1)}n^{f(k+j)}, \quad  0 \le k \le n.
\end{equation} 
 Thus, since all the masses in the inner product are positive, from Theorem \ref{limitecoeficientesKoekoek}, we have 
\begin{align*}\nonumber
&\lim_n \,\lambda_{j,n}=0,\quad j=0, 1,  \dots, r, \\ \nonumber
&\lim_n \,\lambda_{r+1,n}=\lim_n\frac{q_{n}^{(r+1)}(c)}{p_n^{(r+1)}(c)}\frac{1}{A_{r+1}(r+1,n)\,(r+1)!}\not=0.\\ \nonumber
\end{align*}

Now, we will prove that $\lim_n \,\lambda_{r+1,n}=1.$ Indeed, from Corollary \ref{cocientederivadas r+1} and taking into account    (\ref{ligadura}),  we obtain
\begin{align*}
\nonumber
&\frac{q_{n}^{(r+1)}(c)}{p_n^{(r+1)}(c)}\cong \prod_{j=0}^r \frac{(r+1-j)b}{(r+1+j)b+2a+1} \\ \nonumber
&= \prod_{j=0}^r \frac{r+1-j}{r+j+\nu+2}=(r+1)! \,\frac{\Gamma(\nu+r+2)} { \Gamma(\nu+2r+3)}.\\ \nonumber
\end{align*}
Then, from  (\ref{expresionA})  and (\ref{ultimo coeficiente})
$$\lim_n \,\lambda_{r+1,n}=\,\frac{\Gamma(\nu+r+2)} { \Gamma(\nu+2r+3)} \, \lim_n \frac{p_n^{(r+1)}(c)}{p_{n-(r+1)}^{[2(r+1)]}(c)}=1.$$
Finally, it is enough to use the connection formula (\ref{connection-formula}) for the  polynomials  $q_n$ to get their  Mehler-Heine asymptotic.  $\quad \Box$

\medskip

\noindent \textbf{Remark 4}  This theorem is also true when the point $c=\sup I$,  if we change  $(c+z^2/b_n)$ by $(c-z^2/b_n)$ and delete  the factor $(-1)^n$  in formulas (\ref{MH pn-2j}) and  (\ref{MH qn}).

\medskip

\noindent \textbf{Remark 5}  It is important to note that under the hypothesis of Theorem \ref{MH-general}  we have that there exists  the connection formula (\ref{connection-formula}) and moreover    all their  coefficients   tend to zero except the last one which tends to one.

\bigskip

Next, we prove that the hypothesis  (\ref{MH pn-2j}) for all $j=0, 1, \dots, r+1$ in the above theorem can be simplified by certain initial conditions involving only $j=0,1$. To do this,  we will use the following well known formulas:  
\begin{equation}\label{relacion j- j+2}
(x-c)p_{n-1}^{[j+2]}(x)=\frac{\gamma_{n-1}^{[j+2]}}{\gamma_n^{[j]}}\left[ p_n^{[j]}(x)-\frac{p_n^{[j]}(c)}{K_{n-1}^{[j]}(c,c)}\,K_{n-1}^{[j]}(x,c)\right], 
\end{equation}  
\begin{equation}\label{cociente normas j-j+2}
\left (\frac{\gamma_n^{[j]}}{\gamma_{n-1}^{[j+2]}}\right )^2=1+\frac{(p_n^{[j]}(c))^2}{K_{n-1}^{[j]}(c,c)},
\end{equation}
(see  for instance, \cite{gprv})
and 
\begin{equation}\label{relacion nucleo pn-j}
p_n^{[j+1]}(x)=\frac{\gamma_{n}^{[j+1]}}{\gamma_n^{[j]}}\,\frac{K_n^{[j]}(x,c)}{p_n^{[j]}(c)},
\end{equation}
(see \cite{Ch}), where  $\{K_n^{[j]}(x,y) \}_{n}$ denotes the sequence of kernels relative to $\mu_j$.

\bigskip

\begin{proposition} \label{MH 0-1}
  Assume that  the sequence $\{p_n^{[j]}(x) \}_{n\ge0}$ satisfies the asymptotic formulas :
\begin{equation}\label{MH pn-j-0-1}
\lim_{n}\, (-1)^n\frac{a_n^{1/2}}{b_n^{j/2}}\,p_n^{[j]}(c+\frac{z^2}{b_n})=z^{-(\nu+j)}\,J_{\nu+j}(2z), \quad j=0,1
\end{equation}
uniformly  on  compact  subsets of $\mathbb{C},$  where  (\ref {estimacion an-bn}), (\ref{ligadura})  
and  the  conditions
\begin{equation}\label{cociente normas 0-1}
\lim_{n}\, \,\frac{\gamma_n}{\gamma_n^{[1]}}\,\frac{b_n^{1/2}}{n}=\frac{1}{b}
\end{equation}
\begin{equation}\label{cociente normas 1-0}
\lim_{n}\, \,\frac{\gamma_n^{[1]}}{\gamma_{n+1}}\,\frac{b_n^{1/2}}{n}=\frac{1}{b}
\end{equation}
hold.
Then the sequence $\{p_n^{[j]}(x) \}_{n\ge0}$, satisfies  the Mehler-Heine type formula (\ref{MH pn-j-0-1})  for all $j$.
\end{proposition}

\textbf{Proof.}  The idea of the proof is to apply a recursive process such that, whenever we have a Mehler-Heine  type formula  for two consecutive indices $j$ and $j+1$, and moreover the following condition
\begin{equation}\label{cociente normas j/j+1} 
\lim_n \,\frac{\gamma_n^{[j]}}{\gamma_n^{[j+1]}}\,\frac{b_n^{1/2}}{n}=\frac{1}{b}
\end{equation}
holds, then we get  a Mehler-Heine  type formula for $j+2$.

Indeed, suppose that we have (\ref{MH pn-j-0-1}) for $j$ and $j+1$. Then, it can be derived that
\begin{equation*}
(p_n^{[j]})(c)\cong (-1)^{n}\,\frac{A\,B^{j/2}}{\Gamma(\nu+j+1)}n^{f(j/2)}, 
\end{equation*} 
where $f(x)=b\,x+a.$ 
Hence we have
\begin{equation}\label{cociente pn}
 \lim_n \, \frac{p_n^{[j]}(c)}{p_{n-1}^{[j]}(c)} =-1,
 \end{equation}
 and,  since $2f(j/2)+1=b(\nu+j+1)$ by   Stolz criterion,   we obtain
 \begin{equation}\label{cociente pn y nucleo}
K_{n}^{[j]}(c,c)\cong\frac{n(p_n^{[j]}(c))^2}{b(\nu+j+1)}
 \end{equation}
and  from  (\ref{cociente normas j-j+2}) we get
\begin{equation}\label{cociente normas j/j+2}
\gamma_n^{[j]} \cong \gamma_{n-1}^{[j+2]}.
\end{equation}
Now, taking into account  (\ref{relacion j- j+2})  evaluated at  $c+z^2/b_n$ and  (\ref{relacion nucleo pn-j}), we get

\begin{align}\label{auxiliar}
&  \frac{(-1)^{n-1} a_{n}^{1/2}}{b_n^{(j+2)/2}}\, z^2 \,p_{n-1}^{[j+2]}\left(c+\frac{z^2}{b_{n}}\right) =
 \frac{\gamma_{n-1}^{[j+2]}}{\gamma_n^{[j]}} \,\frac{(-1)^{n-1} a_{n}^{1/2}}{b_{n}^{j/2} } \nonumber \\
& \times
\left[ p_{n}^{[j]}\left(c+\frac{z^2}{b_{n}}\right) - \frac{p_{n}^{[j]}(c)\,p_{n-1}^{[j]}(c)}{K_{n-1}^{[j]}(c,c)}\,\frac{\gamma_{n-1}^{[j]}}{\gamma_{n-1}^{[j+1]}}\, p_{n-1}^{[j+1]}\left(c+\frac{z^2}{b_{n}}\right) \right],
\end{align}
and besides, from (\ref{cociente normas j/j+1}), (\ref{cociente pn}) and  (\ref{cociente pn y nucleo}) we have
\begin{align*}
& -\frac{p_{n}^{[j]}(c)\,p_{n-1}^{[j]}(c)}{K_{n-1}^{[j]}(c,c)}\,\frac{\gamma_{n-1}^{[j]}}{\gamma_{n-1}^{[j+1]}}\,b_n^{1/2}  \nonumber \\
&\cong b(\nu+1+j)\,\frac{\gamma_{n-1}^{[j]}}{\gamma_{n-1}^{[j+1]}}\,\frac{b_n^{1/2}}{n} \longrightarrow \nu+1+j.
\end{align*}
To get the Mehler-Heine  type formula for $j+2$,  we only need to take limits in  (\ref{auxiliar}) and  use 
the relation satisfied  by the Bessel functions (see \cite{sz})
$$J_{\nu-1}(z) + J_{\nu+1}(z)  = 2\,\nu\,z^{-1}\,J_{\nu}(z).$$ 
To conclude the proof, it remains to observe that the hypothesis of the proposition is enough to carry out the whole process.
Indeed, the condition (\ref{cociente normas j/j+2})   obtained in each step leads to
$$\gamma_n^{[2j]}\cong \gamma_{n+j},  \quad \gamma_n^{[2j+1]}\cong \gamma_{n+j}^{[1]} .$$ 
So,  from conditions (\ref{cociente normas 0-1}) and (\ref{cociente normas 1-0}), we have $$\lim_{n}\, \,\frac{\gamma_n^{[2j]}}{\gamma_n^{[2j+1]}}  \frac{b_n^{1/2}}{n}=\lim_{n}\, \,\frac{\gamma_{n+j}}{\gamma_{n+j}^{[1]}}\frac{b_n^{1/2}}{n}=\frac{1}{b},$$
and 
$$\lim_{n}\, \,\frac{\gamma_n^{[2j+1]}}{\gamma_n^{[2j+2]}}\frac{b_n^{1/2}}{n}=\lim_{n}\, \,\frac{\gamma_{n+j}^{[1]}}{\gamma_{n+j+1}}\frac{b_n^{1/2}}{n}=\frac{1}{b}.$$
Therefore the required condition (\ref{cociente normas j/j+1}) is satisfied
and  we conclude the proof. $\quad \Box$

\bigskip

\subsection {Symmetric case} 

 Let  $\{p_n(x) \}_{n\ge0}$ be the sequence of orthonormal polynomials with respect to a symmetric measure $\mu$ and $\{ q_n\}_{n\ge0}$ the sequence of orthonormal polynomials with respect to  the inner product  (\ref{innerproduct-simetrico})  where all the masses are positive. Now,
 we will show Mehler-Heine asymptotics for  $\{ q_n\}_{n\ge0}$, using again  a symmetrization process.

As we can see in the following Lemma, there is a remarkable difference with the previous case. We only need to have  the Mehler-Heine formula for $\{ p_n\}$ and an additional condition on the leading coefficients of  $\{ p_n\}$ to achieve a Mehler-Heine formula for $\{p_{n}^{[2j]}\}$ for all $j$.

\begin{lemma} \label{MH-simetrico-2j}  Suppose that the sequence $\{p_n (x) \}_{n\ge0}$ satisfies the Mehler-Heine type formulas :
\begin{align}\label{MH p2n-iniciales}
& \lim_{n}\, (-1)^n \, a_n^{1/2} \,p_{2n} (\frac{z}{b_n})=z^{-\,\nu }\,J_{\nu}(2z)  \\ \nonumber
& \lim_{n}\, (-1)^n\ a_n^{1/2} \,p_{2n+1}(\frac{z}{b_n})=z^{-\,\nu }\,J_{\nu+1}(2z)
 \end{align}
uniformly  on  compact  subsets of $\mathbb{C},$  where 
$$a_n^{-1/2}\cong An^{a}, \quad b_n\cong Bn^b,  \quad A, \,B,\, b>0, \quad  \nu >-1$$
 with 
\begin{equation}\label{ligadura simetrico}
2a+1=2b (\nu+1)
\end{equation}
   and moreover the following conditions hold
\begin{equation}\label{cociente normas pares} 
\lim_n \,\frac{\gamma_{2n}}{\gamma_{2n+1}}\,\frac{b_n}{n}=\frac{1}{2b}
\qquad \text{and} \qquad
\lim_n \,\frac{\gamma_{2n+1}}{\gamma_{2n+2}}\,\frac{b_n}{n}=\frac{1}{2b}.
\end{equation}
Then,  for all $j$,
\begin{align*}
& \lim_{n}\, (-1)^n\frac{a_n^{1/2}}{b_n^{j}}\,p_{2n}^{[2j]}(\frac{z}{b_n})=z^{-\,(\nu+j)}\,J_{\nu+j}(2z) \\ \nonumber
& \lim_{n}\, (-1)^n\frac{a_n^{1/2}}{b_n^j}\,p_{2n+1}^{[2j]}(\frac{z}{b_n})=z^{-\,(\nu+j)}\,J_{\nu+1+j}(2z) 
\end{align*}
uniformly  on  compact  subsets of $\mathbb{C}$.
\end{lemma}

\textbf{Proof.} We will use again the symmetrization process.  From the hypothesis and the relation given in (\ref{simetrizacion pn})
we have
\begin{align*}
& \lim_{n}\, (-1)^n  a_n^{1/2} \,u_{n} (\frac{z^2}{b_n^2})=z^{-\,\nu }\,J_{\nu}(2z)  \\  \nonumber
& \lim_{n}\,  (-1)^n  \frac{ a_n^{1/2} }{(b_n^2)^{1/2}}\,u_n^{[1]}(\frac{z^2}{b_n^2})=z^{-\,(\nu+1)}\,J_{\nu+1}(2z). 
\end{align*}

If we write $u_n^{[j]}(x)= \overline\gamma_n^{[j]}\, x^n+\dots$ the hypothesis (\ref{cociente normas pares}) reads as
\begin{equation*}
\lim_n \,\frac{\overline\gamma_{n}}{\overline \gamma_{n}^{[1]}}\,\frac{(b_n^2)^{1/2}}{n}=\frac{1}{2b}
\qquad \text{and}  \qquad 
\lim_n \,\frac{\overline \gamma_{n}^{[1]}}{\overline\gamma_{n+1}}\,\frac{(b_n^2)^{1/2}}{n}=\frac{1}{2b}.
\end{equation*}

Besides, observe that the condition (\ref{ligadura simetrico}) is now the appropriate to apply Proposition \ref{MH 0-1} to  $\{u_n(x) \}_{n\ge0}$. Therefore, we obtain for all $j$
\begin{equation*}
 \lim_{n}\,  (-1)^n  \frac{a_n^{1/2} }{(b_n^2)^{j/2}}\,u_n^{[j]}(\frac{z^2}{b_n^2})=z^{-\,(\nu +j)}\,J_{\nu+j}(2z),
\end{equation*}
uniformly  on  compact  subsets of $\mathbb{C}.$

Thus, the result follows from (\ref{simetrizacion pn-j}).
 $\quad \Box$

\begin{theorem} \label{MH-simetrico-qn}  With the hypothesis of Lemma  \ref{MH-simetrico-2j},  we have the following Mehler-Heine type formulas for $\{q_n(x) \}_{n\ge0}$
 \begin{align}\label{MH qn-simetrico}
& \lim_{n}\, (-1)^n\,a_n^{1/2}\,q_{2n}(\frac{z}{b_n})= (-1)^{r+1}\, z^{-\,\nu}\,J_{\nu+2r+2}(2z) \\ \nonumber
& \lim_{n}\, (-1)^n\,a_n^{1/2}\,q_{2n+1}(\frac{z}{b_n})= (-1)^{r+1}\,  z^{-\,\nu}\,J_{\nu+2r+3}(2z)
\end{align}
both  uniformly  on  compact  subsets of $\mathbb{C}.$
\end{theorem}

\textbf{Proof.}
From Lemma \ref{MH-simetrico-2j}  and using the following relations described in the symmetrization process
 \begin{equation*}
p_{2n}^{[4j]}(x)=u_n^{[2j]}(x^2), \quad p_{2n+1}^{[4j]}(x)=x\,(u_n^*)^{[2j]}(x^2),
\end{equation*}
we get
\begin{equation*}
 \lim_{n}\, (-1)^n \frac{ a_n^{1/2} }{(b_n^2)^j}\,u_n^{[2j]}(\frac{z^2}{b_n^2})=z^{-\,(\nu +2j)}\,J_{\nu+2j}(2z). 
\end{equation*}
and \begin{equation*}
 \lim_{n}\, (-1)^n  \frac{a_n^{1/2} /b_n}{(b_n^2)^j}\,(u_n^*)^{[2j]}(\frac{z^2}{b_n^2})=z^{-\,(\nu +1+2j)}\,J_{\nu+1+2j}(2z). 
\end{equation*}
To conclude the proof it is enough to check that we can apply Theorem  \ref{MH-general}  to the sequences  $\{u_n(x) \}_{n\ge0}$ and  $\{u_n^*(x) \}_{n\ge0}$.  In this way, we deduce a Mehler-Heine formulas for  $s_n$ and $s_n^*$ and therefore, using (\ref{simetrizacion qn}) for $q_{2n}$ and $q_{2n+1}$, respectively. $\quad \Box$

\medskip

\noindent \textbf{Remark 6}   Note that under the hypothesis of Theorem \ref{MH-simetrico-qn}  we have that there exist  the two connection formulas (\ref{connection-formula-pares}) and  (\ref{connection-formula-impares})  where all their  coefficients tend to zero except the last one in each case that  tends to one.

\subsection  {Asymptotic zero distribution} 

The results above allow us to deduce some additional information about the asymptotic zero distribution  of $\{q_n(x) \}_{n\ge0}$ in terms of the zeros  of the known special functions,  more precisely the  Bessel functions.

Let  $\{p_n(x) \}_{n\ge0}$ be the sequence of polynomials orthonormal with respect to the measure $\mu$ supported on an interval $I$. 
Assume they satisfy, in a neighborhood of the point $c=\inf I$, the following Mehler-Heine formula  :
$$ \lim_{n}\, (-1)^n\,a_n^{1/2}\,p_n(c+\frac{z^2}{b_n})=z^{-\,\nu}\,J_{\nu}(2z),\quad \nu>-1,$$ 
uniformly  on  compact  subsets of $\mathbb{C}.$  This  asymptotic behavior, by Hurwitz's theorem, gives  an additional information of the zeros of $\{p_n(x) \}_{n\ge0}$.  More precisely, if we denote by $ x_{k,n}, \,  1 \le k \le n$  the zeros of the  polynomial $p_n(x)$  in increasing order, and taking into account that the entire function  $z^{-\,\nu}\,J_{\nu}(2z)$ does not vanish at the origin, we can deduce that   for all $k$  
$$  \lim_n \, b_n(x_{k,n}-c) = j_k^{\nu},$$
where  $j_k^{\nu}$  is the $k$th positive zero of $J_{\nu}$. Similar result can be obtained  if we take  $c=\sup I$ where we rename the zeros in decreasing order.

Concerning the zeros of the discrete Sobolev orthogonal polynomials $q_n$, namely $\xi_{k,n}, \,  1 \le k \le n$,  we know that all of them are real and simple and at least $n-(r+1)$ are in the interior of the interval $I$.  Although in a similar way as before,  we have that $\xi_{k,n} \to c$ , there is a remarkable difference in the convergence acceleration of the zeros to $c$. Indeed, since the function limit in  (\ref{MH qn}) has a zero in the origin of multiplicity $r+1$, we have
\begin{align*}
& \lim_n \,b_n(\xi_{k,n}-c) = 0, \qquad  \qquad  1 \le k \le r+1, \\ \nonumber
& \lim_n \, b_n(\xi_{k,n}-c) = j_k^{\nu+2r+2},   \quad  k \ge r+2.
\end{align*}

In the case that the measure $\mu$ is symmetric, since the zeros of the polynomials $q_n$ are symmetric, then analogous information about convergence acceleration of the zeros to $c=0$  can be obtained from formulas (\ref{MH qn-simetrico}).

\bigskip

\section {Examples} 

Let  $\{p_n(x) \}_{n\ge0}$ and $\{q_n(x) \}_{n\ge0}$ be the sequences of orthonormal polynomials with respect to the measure $d\mu(x)$  and  the discrete Sobolev inner product  (\ref{innerproduct})  (or (\ref{innerproduct-simetrico}) if  $\mu$ is symmetric) with all the masses positive.

As  we have seen in  Section 3, under the hypothesis of  Theorems \ref{MH-general}  or \ref{MH-simetrico-qn}  for $\mu$ symmetric, we can assure that there exist the connection formulas  (\ref{connection-formula}) or (\ref{connection-formula-pares})-- (\ref{connection-formula-impares}) where all the coefficients tend to zero except the last one which tends to one. Moreover we have the Mehler-Heine asymptotic  (\ref{MH qn}) or (\ref{MH qn-simetrico}), respectively.

Here, we show some examples of discrete Sobolev polynomials for which the hypothesis of Theorems  \ref{MH-general}  or \ref{MH-simetrico-qn} hold and therefore we get their Mehler-Heine asymptotic in a neighborhood of the point $c$. 

\bigskip

{\bf Laguerre weight}

Let $d\mu(x)=x^{\alpha}\,e^{-x}\,dx, \,\alpha>-1$,  $I=[0, \infty)$,  $c=0$ and  $p_n(x)=l_n^{\alpha}(x)$  the Laguerre orthonormal polynomials.   Since $p_n^{[2j]}(x)=l_n^{\alpha+2j}(x)$, formula (\ref{MH pn-2j}) holds with 
$$a_n^{-1/2}=n^{\alpha/2}, \quad b_n=n,\quad  \nu=\alpha$$ (see \cite{sz}), and therefore we get  (\ref{MH qn}).

\medskip

{\bf Nevai's class} 

Let $ \mu $ be in the well known Nevai's class, $M(0,1)$ and $p_n(x)$ the orthonormal polynomials with respect to $\mu$.  We would like to remark (see, for example, \cite[Theorem 10]{MNT})  that every measure $\mu$ with $\supp \mu=[-1,1]$ and  $\mu^{\prime}>0$ a.e.   ($\mu^{\prime}$ denotes the absolutely continuous part of $\mu$) belongs to $M(0,1)$.

In  \cite[Theorem 1]{A},  for $\mu \in M(0,1)$,  a sufficient condition is given to obtain a Mehler-Heine asymptotic for $p_n(x)$ in a neighborhood of the point  $1$ involving Bessel functions. Notice that,   $\mu_{2j}$ also belongs  to  $M(0,1)$ and so, if the polynomials $p_n^{[2j]}(x)$ satisfy, for all $j=0, 1, \dots, r+1$,  the following condition given in \cite{A}  for $ \nu>-1$
\begin{equation*}
\frac{p_{n+1}^{[2j]}(1)}{p_n^{[2j]}(1)}\cong 1+\frac{\nu+2j+1/2}{n}+o(1/n), \quad  n \to \infty,
\end{equation*} 
then we have that formula (\ref{MH pn-2j}) holds with
\begin{equation}\label{an-bn-Aptekarev}
a_n^{-1/2}=2^{-\nu}\,n^{\nu+1/2}, \quad b_n=n^2/2.
\end{equation}
Thus,  by Theorem  \ref{MH-general} we obtain  the Mehler-Heine asymptotic formula similar to (\ref{MH qn}) with $c=1$. 
Similar results can be obtained for $c = -1$.

\medskip

(a) {\sl Modified Jacobi weight}

Let $ d\mu(x)=h(x) (1-x)^{\alpha}(1+x)^{\beta}\,dx, \, \alpha, \beta>-1$, $I=[-1,1]$  and $h(x)$ a real analytic and positive function on $I$.

Recently in  \cite{F} the author,  using  \cite[Theorem 1]{A},  gives a Mehler-Heine asymptotics for orthonormal polynomials $p_n(x)=p_n^{\alpha, \beta}(x)$  with respect to  $d\mu(x)$  with the restrictions $\alpha >0$ and  $ \beta>0$. 

Thus, for $c=1$,  since $p_n^{[2j]}(x)=p_n^{\alpha+2j, \beta}(x)$ formula (\ref{MH pn-2j}) holds  with $\nu=\alpha$   and the  values of $a_n$ and $b_n$ are given by (\ref{an-bn-Aptekarev}). Then,  we get the analogous formula to (\ref{MH qn})  with the corresponding changes due to $c=\sup I$, see Remark 4.

On the other hand, for $c=-1$,  since $p_n^{[2j]}(x)=p_n^{\alpha, \beta+2j}(x)$  formula (\ref{MH pn-2j}) holds with the same values of $a_n$ and $b_n$ but now with $\nu=\beta$. Therefore, we can assure that  (\ref{MH qn}) holds.

\medskip

(b) {\sl Jacobi weight}

Let $ d\mu(x)= 2^{\alpha- \beta} (1-x)^{\alpha}(1+x)^{\beta}\,dx, \, \alpha>-1, \, \beta>-1$, $I=[-1,1]$.   The same asymptotics are true for $c=\pm 1 $  but now there are no restrictions on the parameters $\alpha$ and $\beta$, that is $\alpha, \beta >-1$ (see \cite{sz}).

\medskip
  
{\bf Freud weight}

Let $d\mu(x)=\exp(-2\vert x \vert^{\alpha})\,dx, \,\alpha>1$,  $I=\mathbb{R}$, $c=0$ and $p_n(x)$ the Freud orthonormal polynomials with respect to $d\mu(x)$.  A  Mehler-Heine asymptotic near to $0$ for $p_n$ was obtained in \cite[Theorem 1]{ampr}. From there,  with a change of notation, we have  that  (\ref{MH p2n-iniciales}) holds with $\nu=-1/2$ and
 $$a_n^{-1/2}=\sqrt{2}\, (2c_{\alpha})^{-1/2{\alpha}}\,n^{-1/2{\alpha}}, \quad  b_n = \frac{\alpha}{\alpha-1}(2c_{\alpha})^{-1/{\alpha}}\,n^{1-1/{\alpha}},$$
where $c_{\alpha}=\frac{\sqrt{\pi}}{2}\frac{\Gamma(\alpha/2)}{\Gamma((\alpha+1)/2)}$.

Moreover, from the result given in \cite[p. 365]{ST97} it is easy to deduce that the conditions (\ref{cociente normas pares}) hold. Then, since  condition (\ref{ligadura simetrico}) is trivially satisfied,  we can apply Theorem \ref{MH-simetrico-qn}, and therefore we get  (\ref{MH qn-simetrico}).


\bigskip


 
\begin{thebibliography}{10}

\bibitem{ammr}
M. Alfaro, F. Marcell\'{a}n, H.G. Meijer, M.L. Rezola, {\em Symmetric orthogonal polynomials for Sobolev-type inner products}, 
\newblock J. Math. Anal. Appl.  {\bf 184 } (1994),  360-381.

\bibitem{ampr2011}
M. Alfaro, J.J. Moreno-Balc\'{a}zar, A. Pe\~{n}a, M.L. Rezola, {\em A new approach to the asymptotics of Sobolev type orthogonal polynomials}, 
\newblock J.  Approx. Theory {\bf 163} (2011), 460-480.

\bibitem{ampr}
M. Alfaro, J.J. Moreno--Balc\'{a}zar, A.  Pe\~{n}a, M.L. Rezola, {\em Asymptotic formulae for generalized Freud polynomials}, 
\newblock J. Math. Anal. Appl.   {\bf 421 } (2015),  474-488.

\bibitem{am}
R. \'{A}lvarez-Nodarse, J.J. Moreno--Balc\'{a}zar,  {\em Asymptotic properties of  generalized Laguerre orthogonal  polynomials}, 
\newblock Indag. Math. (N.S.)  {\bf 15 } (2004),  151-165.


\bibitem{A}
A.I. Aptekarev, {\em Asymptotics of orthogonal polynomials in a neighborhood of  the endpoints of the interval of orthogonality}, 
\newblock Russian Acad. Sci. Sb. Math. {\bf 76 } (1993),  35-50.

\bibitem{Ch} 
T.S. Chihara, {\em An Introduction to Orthogonal Polynomials},  \newblock Gordon and Breach, New York, 1978.


\bibitem{F}
B. Xh. Fejzullahu, {\em Mehler-Heine formulas for orthogonal polynomials with respect to the modified Jacobi weight},
\newblock Proc. Amer. Math. Soc.  {\bf 142} (2014), 2035--2045.

\bibitem{gprv}  J.J. Guadalupe, M. P\'erez, F.J. Ruiz, J.L. Varona, {\em Asymptotic behaviour of orthogonal polynomials relative to measures with mass points},
 \newblock  Mathematika  {\bf 40} (1993),  331--344.

\bibitem{k}
G. Klambauer,  {\em Aspects of Calculus}, 
\newblock  Springer--Verlag, New York, 1986.

\bibitem{k1}
R. Koekoek, {\em Generalizations of Laguerre polynomials}, 
\newblock J. Math. Anal. Appl.  {\bf 153} (1990), 576--590.

\bibitem{km}
R. Koekoek, H.G. Meijer, {\em A generalization of Laguerre polynomials}, 
\newblock SIAM J. Math. Anal.  {\bf 24}  (1993),  768--782.

\bibitem{K-Mc}  T. Kriecherbauer, K.~T.--R. McLaughlin,
{\em Strong asymptotics of  polynomials orthogonal with respect to Freud weights},
\newblock Internat. Math. Res. Notices,  {\bf 6}  (1999), 299--333.

\bibitem{LMV}  G. L\'{o}pez,  F. Marcell\'{a}n, W. Van Assche, {\em Relative asymptotics for polynomials orthogonal with respect to a discrete Sobolev inner product},
 \newblock  Constr. Approx. {\bf 11} (1995),  107--137.



\bibitem{mzfh}
F. Marcell\'{a}n, R. Zejnullahu, B. Fejzullahu,  E. Huertas  {\em On orthogonal polynomials with respect to certain discrete Sobolev inner product}, 
\newblock Pacific J. Math.  {\bf 257} (2012), 167--188.


\bibitem{MNT}  A. M\'{a}t\'{e}, P.  Nevai,  V. Totik, {\em Extensions of Szeg\H{o}'s Theory of  Orthogonal Polynomials II},
 \newblock  Constr. Approx. {\bf 3} (1987),  51--72.

\bibitem{pr}
A.  Pe\~{n}a, M.L. Rezola, {\em Discrete Laguerre-- Sobolev expansions: A Cohen type inequality}, 
\newblock J. Math. Anal. Appl.  {\bf 385} (2012),  254--263.

\bibitem{rms}
I.A. Rocha, F. Marcell\'{a}n, L. Salto,   {\em  Relative asymptotics and  Fourier series of orthogonal polynomials with a discrete Sobolev  inner product}, 
\newblock J.  Approx. Theory {\bf 121} (2003), 336--356.

\bibitem{ST97}
E.B. Saff, V.  Totik,
 {\em Logarithmic potentials with external fields},
 Grundlehren der Mathematischen Wissenschaften, {\bf 316}.
Springer-Verlag, Berlin, 1997.

\bibitem{sz}
G.~Szeg\H{o}, {\em  Orthogonal Polynomials}, \newblock Amer. Math. Soc.
Colloq. Publ. vol.  23, Amer. Math. Soc., Providence RI,
1975.

\bibitem{W}
W.  Van Assche, {\em Mehler-Heine asymptotics for multiple orthogonal polynomials}, 
\newblock arXiv:1408.6140v1(2014).

\end{thebibliography}
\end{document}